\documentclass[10pt,a4paper]{article}

\usepackage{amssymb}
\usepackage{amsmath}
\usepackage{amsthm}
\usepackage{tikz}
\usepackage{accents}
\usetikzlibrary{shapes}
\usepackage{textcomp}
\usepackage{caption}
\usepackage{undertilde}
\usepackage{ulem}
\newcommand{\vect}[1]{\boldsymbol{#1}} 
\newcommand{\tensor}[1]{\boldsymbol{#1}} 

\begin{document}

\title{\bfseries{Coupling of complex function theory and finite element method for crack propagation through energetic formulation: conformal mapping approach and reduction to a Riemann-Hilbert problem}}

\author{Dmitrii Legatiuk\thanks{Chair of Applied Mathematics, Bauhaus-Universit\"at Weimar, Coudraystr. 13B, 99423 Weimar, Germany} \and Daniel Weisz-Patrault\thanks{LMS, \'{E}cole Polytechnique, Institut Polytechnique de Paris, F-91128 Palaiseau, France}}
\date{}

\maketitle

\begin{abstract}
In this paper we present a theoretical background of a coupled analytical-numerical approach to model a crack propagation process in two-dimensional bounded domains. The goal of the coupled analytical-numerical approach is to obtain the correct solution behaviour near the crack tip by help of the analytical solution constructed by using tools of the complex function theory and couple it continuously with the finite element solution in the region far from singularity. In this way, crack propagation could be modelled without using remeshing. Possible directions of crack growth can be calculated through the minimization of the total energy composed of the potential energy and the dissipated energy based on the energy release rate. Within this setting, an analytical solution of a mixed boundary value problem based on complex analysis and conformal mapping techniques is presented in a circular region containing an arbitrary crack path. More precisely, the linear elastic problem is transformed into a Riemann-Hilbert problem in the unit disk for holomorphic functions. Utilising advantages of the analytical solution in the region near the crack tip, the total energy could be evaluated within short computation times for various crack kink angles and lengths leading to a potentially efficient way of computing the minimization procedure. To this end, the paper presents a general strategy of the new coupled approach for crack propagation modelling. Additionally, we also discuss obstacles on the way of practical realisation of this strategy.
\end{abstract}

\section{Introduction}

Methods of complex function theory provide various tools to construct exact solutions to differential equations, especially in the case of singularity, such as e.g. crack tip problem in linear elastic fracture mechanics. Particularly, with the introduction of famous Kolosov-Muskhelishvili formulae, methods of complex function theory became indispensable to handle problems of linear elasticity \cite{Muskhelishvili_2}. The classical Kolosov-Muskhelishvili formulae enable us to represent displacements and stresses of a two-dimensional elastic body in terms of two holomorphic functions $\Phi(z)$ and $\Psi(z)$, $z\in\mathbb{C}$. Because of obvious advantages of the function-theoretic approach, such as exact singular behaviour near the crack tip and preservation of all basic physical assumptions, methods of complex function theory constituted the foundation of classical fracture mechanics \cite{Liebowitz,sif}.\par
A known disadvantage of function-theoretic methods is the fact that a complete boundary value problem can be solved explicitly only for some elementary (simple) domains, such as e.g. the unit disk or half-plane. Considering that domain coming from real-world engineering problems have generally more complicated geometry, numerical methods, such as e.g. extended finite element method \cite{Moes99}, are frequently used to solve static and dynamic fracture mechanics problems nowadays. The idea of modern numerical methods used in fracture mechanics application is to enrich classical finite element shape functions with known analytical solution, e.g. Westergaard solution and partition of unity \cite{Moes99}, to obtain correct asymptotic behaviour near the crack tip. The drawback of such methods is the lost continuity between enriched and standard elements, since the modified shape functions do not satisfy the interpolation conditions. Thus, the methods obtained in this way do not satisfy basic assumptions of the classical theory of finite element method \cite{Ciarlet}, and therefore, it is difficult to perform a rigorous convergence analysis.\par
In this context, utilising advantages of both function-theoretic methods and finite element method, coupled analytical-numerical methods could be alternative approaches towards higher accuracy of solutions in the region near the singularity. While a coupling between the analytical solution obtained by function-theoretic methods and the finite element solution can be introduced in several ways (see e.g. \cite{Piltner_1,Piltner_2}), we focus on a continuous coupling in this paper. The idea of a continuous analytical-numerical coupling is to introduce a special interpolation operator preserving $C^{0}$ continuity of the displacement field on the interface between function-theoretic solution and the classical finite elements. Construction of such an interpolation operator has been presented in \cite{Guerlebeck_1,Guerlebeck_2}, and convergence analysis of the coupled method has been performed in \cite{Guerlebeck_3,Legatiuk_1}, where the coupling error has been also estimated explicitly. However, only problems of fracture mechanics with static cracks have been considered so far. Therefore, in this paper, we present an extension of the coupled approach to crack propagation problems in two-dimensional domains.\par
The crack propagation approach presented in this paper is within the framework of linear elastic fracture mechanics. The main result of this theory of fracture is that linear elastic calculations are sufficient to estimate the fracture energy release rate, or equivalently the stress intensity factor, to determine whether a crack propagates or not. However, a prediction of the crack propagation path with bifurcation points cannot be obtained only by considering the fracture energy release rate. Therefore, additional bifurcation criteria have been introduced for computing the crack propagation path, e.g. the maximum hoop stress criterion \cite{erdogan1963}, the maximum-energy-release-rate criterion \cite{MAXG} or an asymptotic expansion of the stress intensity factor \cite{leblond1989crack}. Moreover, one of the most elegant and physically consistent approaches is the variational formulation proposed in \cite{francfort-marigo}. In this approach, at each time, and for any boundary condition, the crack propagation path is obtained by finding the global minimum of the total energy, which is the sum of potential energy and dissipated energy, under the constraint of the irreversibility of the crack growth. One of the main advantages of the energetic formulation is its correspondence to a quasi-static evolution of the crack, implying that a succession of stable states is simulated without referring to the detailed mechanisms arising between two stable states. Thus, the energetic approach authorises discontinuous evolutions, practically meaning that, for instance, the crack length increase is not infinitesimal but may be finite between two time steps. However, the minimisation procedure may be time-consuming as numerical methods (that should be sufficiently refined to obtain acceptable accuracy) are repeatedly used to minimise the total energy.\par
In this paper, an approach combining coupled analytical-numerical method and energetic approach in order to model crack propagation is proposed. The expected advantage of such an approach is a reduced computation time on the finite element side, since the analytical solution near the crack tip is used. However, since the original coupled analytical-numerical method is limited to cracks without bifurcation points \cite{Guerlebeck_3}, the analytical solution at first must be extended to the case of kinked cracks. This extension is done by using a conformal mapping approach, and therefore, the linear elastic problem in the region near the crack tip is reduced to a Riemann-Hilbert problem for holomorphic functions in the unit disk. Therefore, our aim here is to extend the conformal mapping approach to the case of coupled analytical-numerical method. As it will be discussed in the paper, practical (numerical) realisation of this approach still needs to be addressed properly due to known difficulties with numerical conformal mappings. Therefore, this paper aims at presenting a general strategy for modelling crack propagation based on a continuous coupling of function-theoretic methods and the finite element method. Moreover, we present an explicit solution of the Riemann-Hilbert problem and provide a 
detailed discussion on future steps for practical realisation of the proposed method along with first numerical calculations for the Riemann-Hilbert problem.\par

\section{Modelling crack propagation via the coupling of function-theoretic and finite element methods}

In this section we present a general description of the method to model crack propagation via a continuous coupling of complex function theory and finite element method. To support the reader, we start with a general overview of the coupled method underlying only essential steps relevant for the crack propagation modelling. After that, we discuss the mechanical point of view on the propagation process and outline the idea to use conformal mapping approach leading to the formulation of a Riemann-Hilbert problem, which is discussed in details in the upcoming sections.\par

\subsection{Continuous analytical-numerical coupling for static cracks}

Let $G\subset\mathbb{C}$ be a simply connected bounded domain containing a crack. Further, let $\Gamma$ be a boundary of $G$, and it is assumed to be sufficiently smooth except the turning point given by a crack tip, which causes a well-known crack-tip singularity. We consider now the classical boundary value problem of linear elasticity formulated as follows
\begin{equation}
\label{linear_elasticity}
\left\{
\begin{array}{rcll}
-\mu\,\Delta\mathbf{u}-(\lambda+\mu)\,\mathrm{grad}\,\mathrm{div}\,\mathbf{u} & = & \mathbf{f} & \mbox{in } G,\\
\mathbf{u} & = & \mathbf{g}_{0} & \mbox{on } \Gamma_{0},\\
\displaystyle \boldsymbol{\sigma}\cdot \overline{n} & = & \mathbf{g}_{1} & \mbox{on } \Gamma_{1}, 
\end{array}\right.
\end{equation}
where $\lambda$ and $\mu$ are classical Lam\'e constants, and $\mathbf{f}$ is the density of volume forces, $\mathbf{u}$ is the unknown displacement vector, $\boldsymbol{\sigma}$ is the Cauchy stress tensor, $\overline{n}$ is the unit outer normal vector, and $\Gamma_{0}$ and $\Gamma_{1}$ are parts of the boundary with Dirichlet and Neumann boundary conditions ($\mathbf{g}_{0}$ and $\mathbf{g}_{1}$), respectively.\par
\begin{figure}[h!]
\centering
\includegraphics[width=1.0\textwidth]{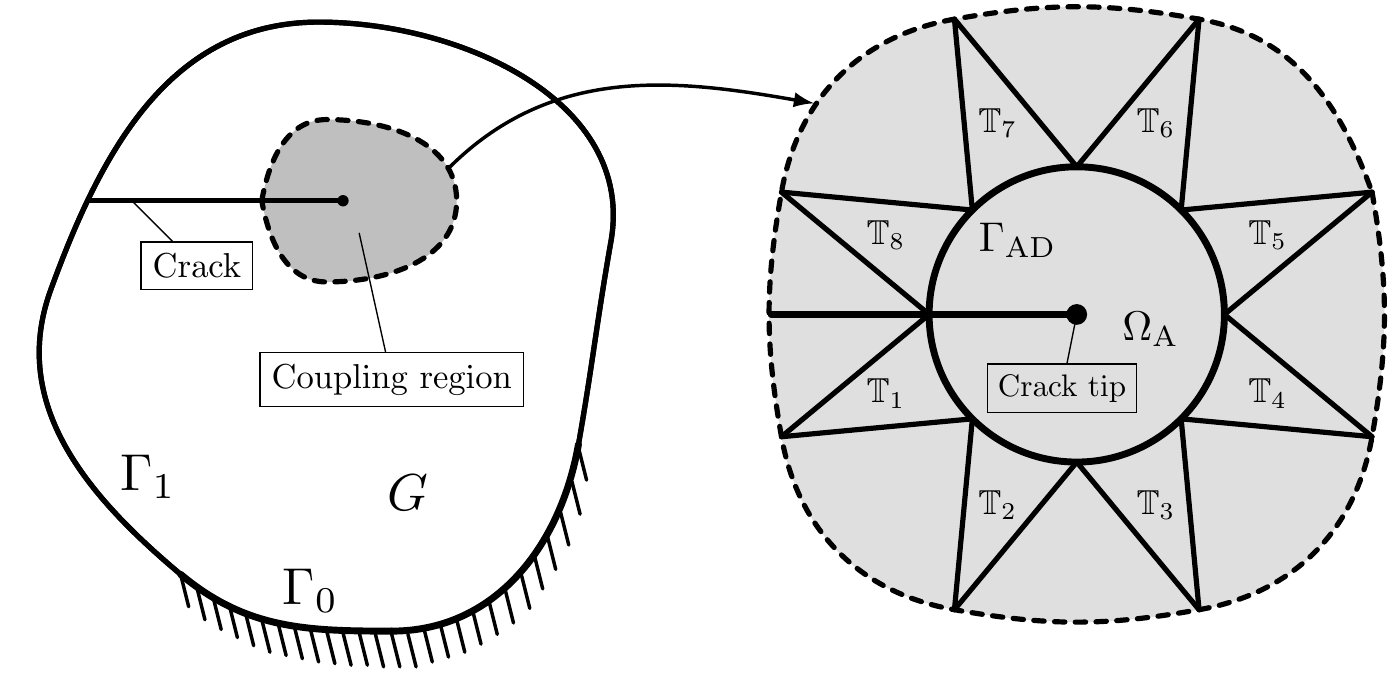}
\caption{Left: domain $G$ containing a crack and the coupling region. Right: further subdivision of the coupling region into analytical domain $\Omega_{\mathrm{A}}$ and coupling elements $\mathbb{T}_{i}, i=1,\ldots,8$.}
\label{figure_1}
\end{figure}

To provide an exact description of the solution behaviour near the singularity, we introduce a local {\itshape coupling region} surrounding the crack tip (see Fig. \ref{figure_1}, left). The right side of Fig. \ref{figure_1} illustrates the coupling region with more details. Particularly, the coupling region is further subdivided into {\itshape analytical domain}  $\Omega_{\mathrm{A}}$ circled by curved triangular elements $\mathbb{T}_{i}$ (8 in Fig. \ref{figure_1}), which are called {\itshape coupling elements}. The interface $\Gamma_{\mathrm{AD}}$ between $\Omega_{\mathrm{A}}$ 
and coupling elements is called the {\itshape coupling interface}. The remaining part of the domain $G$ is triangulated by standard finite elements.\par
Introducing the coupling region enables us to couple continuously the exact solution to the differential equation of linear elasticity in $\Omega_{\mathrm{A}}$ with finite element solution in the remaining part of the domain. This continuous coupling is provided by help of a special interpolation operator, which is based on the analytical solution. A detailed construction of such an interpolation operator and its invariance property have been discussed in \cite{Guerlebeck_1,Guerlebeck_2}. Because of the continuous coupling, a variational problem as in the classical finite element method (FEM) theory, see \cite{Ciarlet} for details, can be formulated in our case. Since the goal of this paper is not to discuss finite element aspects of the coupled method, but rather focus on function-theoretic tools to model crack propagation, we omit all further technical details on the FEM part of the method and refer to \cite{Legatiuk_2} for a complete construction.\par
Analytical solution to the differential equation in $\Omega_{\mathrm{A}}$ is constructed by the Kolosov-Muskhelishvili formulae \cite{Muskhelishvili_2}; in polar coordinates these formulae allow us to represent components of the displacement field and stress tensor in the following form
\begin{equation*}
\begin{array}{rcl}
\displaystyle 2\mu(u_{r}+i\,u_{\varphi}) & = & \displaystyle e^{-i\varphi}\left(\kappa\,\Phi(z)-z\,\overline{\Phi'(z)}-\overline{\Psi(z)}\right), \\
\\
\displaystyle \sigma_{rr}+\sigma_{\varphi\varphi} & = & \displaystyle 2\left[\Phi'(z)+\overline{\Phi'(z)}\right], \\
\\
\displaystyle \sigma_{\varphi\varphi}-\sigma_{rr}+2i\,\sigma_{r\varphi} & 

= & \displaystyle 2e^{2i\,\varphi}\left[\bar{z}\,\Phi''(z)+\Psi'(z)\right],
\end{array}
\end{equation*}
where $\Phi(z)$ and $\Psi(z)$ are two holomorphic functions, and $\kappa\in(1,3)$ is the Kolosov's constant. For static cracks, the holomorphic function $\Phi(z)$ and $\Psi(z)$ have been written in terms of power series expansion \cite{Guerlebeck_1}
\begin{equation*}
\Phi(z)=\sum_{k=0}^{\infty}a_{k}z^{\lambda_{k}},\quad \Psi(z)=\sum_{k=0}^{\infty}b_{k}z^{\lambda_{k}}, \mbox{ with } a_{k},b_{k}\in\mathbb{C}, \lambda_{k}\in\mathbb{R}.
\end{equation*}
Using these series expansions in the Kolosov-Muskhelishvili formulae and applying traction free boundary conditions on the crack faces, exponents $\lambda_{k} = \frac{k}{2}$, $k=1,2,\ldots$ are found, which correspond to the classical crack tip singularity, see \cite{Liebowitz} for details. Moreover, relations between complex coefficients $a_{k}$ and $b_{k}$ are also identified, and therefore, the displacement field can be written now as follows 
\begin{equation}
\label{displacement_cartesian_coordinates_with_coeff}
\begin{array}{rcl}
2\mu(u_{1}+i\,u_{2}) & = & \displaystyle \sum\limits_{n=0,2,...}^{\infty}r^{\frac{n}{2}}\left[a_{n} \left(\kappa\,e^{i\varphi\frac{n}{2}}+e^{-i\varphi\frac{n}{2}}\right) + \right. \\ 
\\
& & \displaystyle + \left. \frac{n}{2}\bar{a}_{n}\left(e^{-i\varphi\frac{n}{2}}-e^{-i\varphi(\frac{n}{2}-2)}\right)\right] + \\
\\
& & \displaystyle + \sum\limits_{n=1,3,...}^{\infty}r^{\frac{n}{2}}\left[a_{n} \left(\kappa\,e^{i\varphi\frac{n}{2}}-e^{-i\varphi\frac{n}{2}}\right) + \right. \\
\\
& & \displaystyle \left. + \frac{n}{2}\bar{a}_{n}\left(e^{-i\varphi\frac{n}{2}}-e^{-i\varphi(\frac{n}{2}-2)}\right)\right],
\end{array}
\end{equation}
where unknown coefficients $a_{n}$ are still to be identified by solving the global boundary value problem~(\ref{linear_elasticity}) via the coupled finite element procedure.\par
Finally, the continuity of displacement field through the entire coupling interface $\Gamma_{\mathrm{AD}}$ in the finite element procedure is preserved by constructing finite element basis functions based on the truncated exact solution~(\ref{displacement_cartesian_coordinates_with_coeff}). Let us consider $n$ nodes on the interface $\Gamma_{\mathrm{AD}}$ belonging to the interval $[-\pi,\pi]$, then the interpolation function $f_{n}(\varphi)$ restricted to $\Gamma_{\mathrm{AD}}$, i.e. $r=r_{\mathrm{A}}$, has the following form
\begin{equation}
\label{interpolation_function_general_form}
\begin{array}{lcl}
f_{n}(\varphi) & = & \displaystyle \sum\limits_{k=0,2,...}^{N_{1}}r_{\mathrm{A}}^{\frac{k}{2}}\left[a_{k} \left(\kappa\,e^{i\varphi\frac{k}{2}}+e^{-i\varphi\frac{k}{2}}\right) +\frac{k}{2}\bar{a}_{k}\left(e^{-i\varphi\frac{k}{2}}-e^{-i\varphi(\frac{k}{2}-2)}\right)\right]+\\
\\
& & +\displaystyle\sum\limits_{k=1,3,...}^{N_{2}}r_{\mathrm{A}}^{\frac{k}{2}}\left[a_{k} \left(\kappa\,e^{i\varphi\frac{k}{2}}-e^{-i\varphi\frac{k}{2}}\right) +\frac{k}{2}\bar{a}_{k}\left(e^{-i\varphi\frac{k}{2}}-e^{-i\varphi(\frac{k}{2}-2)}\right)\right],
\end{array}
\end{equation}
where the numbers of basis functions $N_{1}$ and $N_{2}$ are related to $n$ as follows:
\begin{equation*}
\begin{array}{ll}
N_{1} = \left\{\begin{array}{ll} n-2 & \mbox{for even } n,\\
n-1 & \mbox{for odd } n,\end{array}\right. & 
N_{2} = \left\{\begin{array}{ll} n-1 & \mbox{for even } n,\\
n-2 & \mbox{for odd } n.\end{array}\right.
\end{array}
\end{equation*}
The basis functions for finite element approximation are then obtained by interpolating the unknown displacements $\mathbf{U}_{j}$, $j=0,\ldots,n-1$ on the coupling interface $\Gamma_{\mathrm{AD}}$, see  \cite{Guerlebeck_3,Legatiuk_1,Legatiuk_2} for all further details.\par
In summary, this paper aims at extending this coupling strategy to crack propagation, which implies to consider more complex crack paths and therefore to develop an adapted analytical solution in the analytical domain $\Omega_{\mathrm{A}}$.\par

\subsection{Strategy to model crack propagation}

A typical approach to model crack propagation by help of the finite element method is based on the idea of a local or global remeshing at each step of crack propagation. Although this approach can be immediately adapted to our setting, it is well-known that remeshing is computationally costly and inefficient. Alternatively, we prefer to utilise the advantage of the coupled method enabling us to work with a fixed size of the analytical domain $\Omega_{\mathrm{A}}$ without involving a global refinement. In this case, we allow the crack to propagate only inside the analytical domain that should be taken as large as possible, while performing refinement on the mesh around $\Omega_{\mathrm{A}}$.\par
\begin{figure}[h!]
\centering
\includegraphics[width=1.0\textwidth]{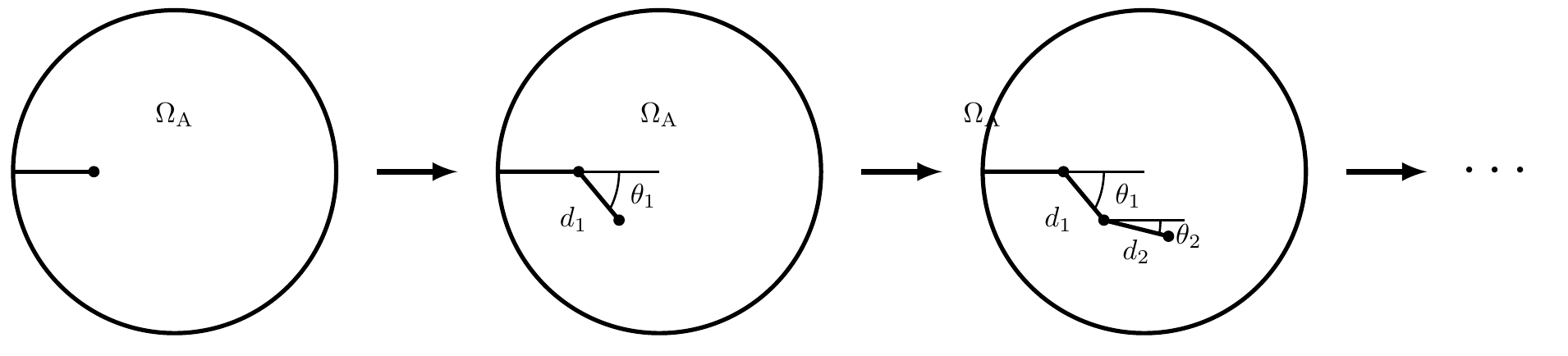}
\caption{Development of the crack inside the analytical domain $\Omega_{\mathrm{A}}$ for first few loading steps}
\label{figure_3_crack_development}
\end{figure}

Let us now consider more precisely the analytical domain $\Omega_{\mathrm{A}}$. At the initial moment, the crack tip is located inside $\Omega_{\mathrm{A}}$, and the crack faces are going along the negative direction of $x_{1}$-axis of a Cartesian coordinate system. After the first loading step, the crack is allowed to propagate inside the analytical domain. We assume that the crack propagates with a finite length at one loading step, i.e. the crack tip moves along the propagation direction defined by the angle $\theta_{i}$ for a finite length $d_{i}$ with $i=1,2,\ldots$ denoting the loading step, see Fig. \ref{figure_3_crack_development}. To evaluate the angle $\theta_{i}$ and the length $d_{i}$ we have to solve a minimisation problem according to \cite{francfort-marigo}, and therefore, to construct an analytical solution to the crack tip problem in $\Omega_{\mathrm{A}}$.\par 
As already mentioned, the analytical solution~(\ref{displacement_cartesian_coordinates_with_coeff}) cannot be used to calculate the displacement field for next loading steps, since the basic assumptions of the model are not satisfied any more due to the presence of a kinked crack. This problem can be solved by application of a conformal mapping, which allows us to map the analytical domain after several loading steps (see Fig.~\ref{figure_3_crack_development}) to the unit disk. The solution of a boundary value problem in the unit disk can be obtained again by the Kolosov-Muskhelishvili formulae. According to \cite{Muskhelishvili_2}, these the Kolosov-Muskhelishvili formulae under a conformal mapping are written as follows
\begin{equation}
\label{kolosov_formulae_under_mapping}
\begin{array}{rcl}
\displaystyle \sigma_{rr}+\sigma_{\varphi\varphi} & = & \displaystyle 2\left[\Phi(\zeta)+\overline{\Phi(\zeta)}\right], \\
\\
\displaystyle \sigma_{rr}+i\,\sigma_{r\varphi} & = & \displaystyle \Phi(\zeta)+\overline{\Phi(\zeta)} - \frac{\bar{\zeta}^{2}}{r^{2}\omega'(\zeta)}\left[\omega(\zeta)\overline{\Phi'(\zeta)} + \overline{\omega'(\zeta)}\,\overline{\Psi(\zeta)}\right],\\
\\
\displaystyle 2\mu|\omega'(\zeta)|(u_{r}+i\,u_{\varphi}) & = & \displaystyle \frac{\bar{\zeta}}{r}\overline{\omega'(\zeta)}\left[\kappa\,\eta(\zeta) - \omega(\zeta)\overline{\Phi(\zeta)}-\overline{\chi(\zeta)}\right],

\end{array}
\end{equation}
where $r,\varphi$ denote polar coordinates in the unit disk, $\zeta=r\exp(i\varphi)$, $\Phi(\zeta)$ and $\Psi(\zeta)$ are two holomorphic functions defined on the unit disk and, and $\eta(\zeta)$ and $\chi(\zeta)$ are functions related to $\Phi(\zeta)$, $\Psi(\zeta)$ by help of the expressions
\begin{equation}
\label{relation_between_functions}
\eta'(\zeta)=\Phi(\zeta)\omega'(\zeta), \qquad \chi'(\zeta)=\Psi(\zeta)\omega'(\zeta),
\end{equation}
and $\omega(\zeta)$ is a mapping from the original geometry to the unit disk. Solution of a boundary value in the unit disk and construction of a mapping $\omega(\zeta)$ is described in detail in Section \ref{section_conformal_mapping}.\par
In addition, if $1,2$ denote Cartesian directions in the original geometry, displacements $u_{1}$, $u_{2}$ read according to \cite{Muskhelishvili_2} as follows
\begin{equation}
\label{displacement_cartesian}
2\mu(u_{1}+i\,u_{2})  = \kappa\,\eta(\zeta) - \omega(\zeta)\overline{\Phi(\zeta)}-\overline{\chi(\zeta)}.
\end{equation}

\section{Conformal mapping for a cracked disk and the Riemann-Hilbert problem}\label{section_conformal_mapping}

In this section, we discuss the application of conformal mapping to construct an analytical solution for a crack disk and formulation of the corresponding Riemann-Hilbert problem in the unit disk. Moreover, to keep construction general, we do not specify the corresponding conformal mapping explicitly, although the classical Schwarz-Christoffel mapping is the first candidate \cite{Driscoll}. We come back to this point later during the discussion in Section~\ref{Section_conclusions}.\par

\subsection{Application of the conformal mapping to a cracked disk}

The idea of using conformal mapping for studying crack propagation within domain $\Omega_{\mathrm{A}}$ is motivated by several facts: (i) analytical solution~(\ref{displacement_cartesian_coordinates_with_coeff}) is not valid for the case of a propagated crack, since the distance between the crack tip and the kinking point is too small to validate the assumptions of the classical crack tip solution; (ii) remeshing is not necessary if propagating crack does not intersect the coupling interface $\Gamma_{\mathrm{AD}}$; (iii) analytical constructions are expected to provide higher flexibility and accuracy in calculating mechanical quantities of interest relevant for propagation process \cite{sif}.\par
\begin{figure}[h!]
\centering
\includegraphics[width=1.0\textwidth]{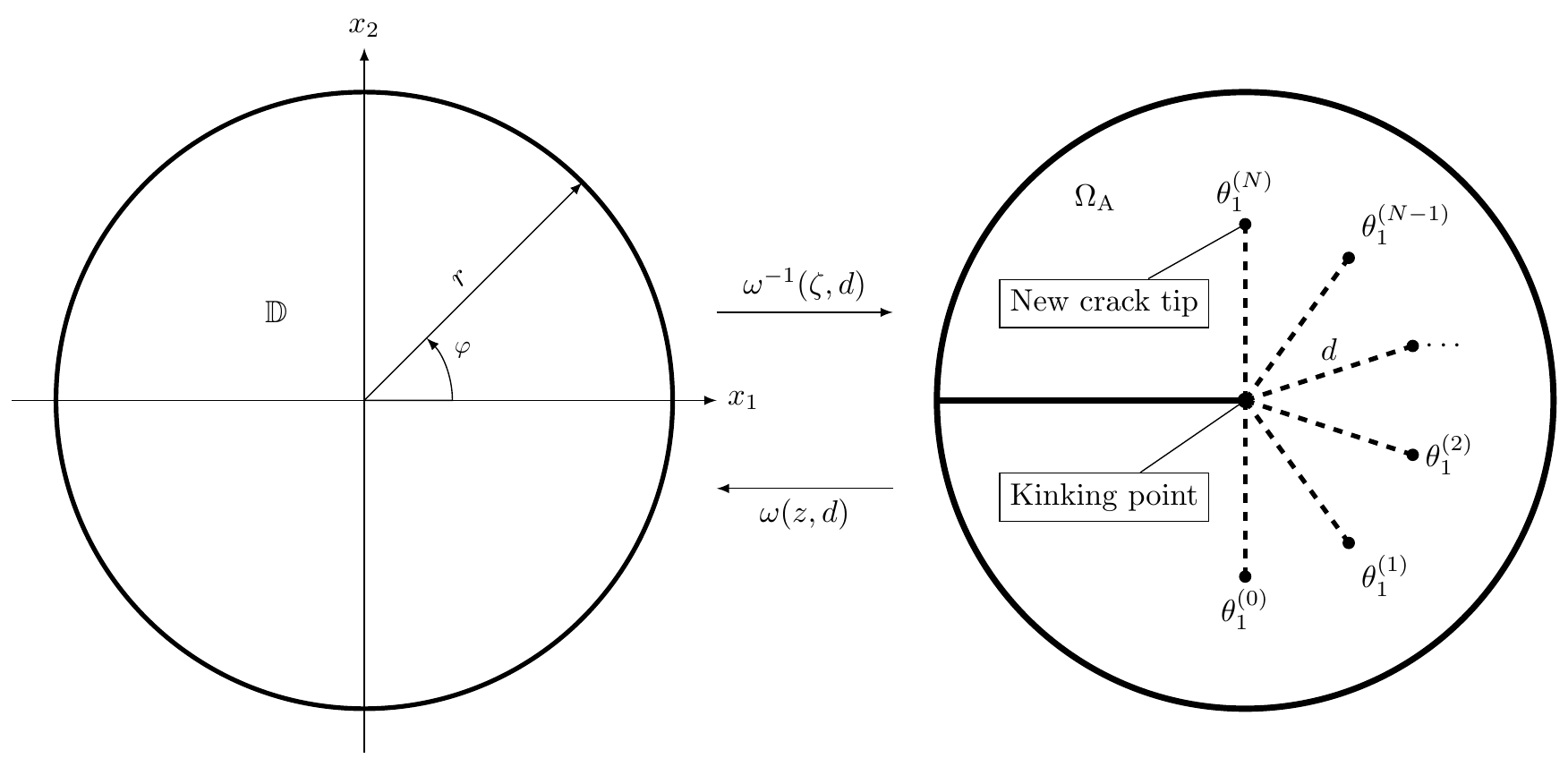}
\caption{Mapping between the unit disk and a cracked disk with different possible directions for crack propagation}
\label{figure_5}
\end{figure}

Looking at the crack propagation process from the mechanical point of view, it is known that depending on specific loading conditions the crack can propagate in different directions controlled by the angle $\theta_{i}$ with the propagation length $d_{i}$, where $i$ is the number of loading step. Practically it implies, that conformal mappings need to be calculated for all possible directions and lengths, which is a computationally expensive operation to perform online. However, considering that the crack propagates only inside $\Omega_{\mathrm{A}}$, conformal mappings can be pre-calculated for different values of the angle $\theta_{i}^{(k)}\in\left[-\frac{\pi}{2},\frac{\pi}{2}\right]$ with $i$ denoting propagation step and $k=0,\ldots,N$ being the number of a specific angle with parameter $N$ controlling the angular discretisation, and having the lengths $d$ as a free parameter in the mapping, see Fig. \ref{figure_5}. Note that the crack tip is located at the centre of $\Omega_{\mathrm{A}}$ in Fig.~\ref{figure_5} only for clarity reasons. In practice, it is better to place the crack tip sufficiently close to the boundary of $\Omega_{\mathrm{A}}$ (taking into that traction-free assumptions on the crack faces must be still satisfied) for addressing more propagation steps inside $\Omega_{\mathrm{A}}$ without remeshing.\par

\subsection{Boundary value problem of linear elasticity as a Riemann-Hilbert problem}

Now we will show how a boundary value problem of elasticity can be transformed into a Riemann-Hilbert problem for a piecewise holomorphic function. Let now $\mathbb{D}=\left\{\zeta \in\mathbb{C}\colon |\zeta|<1\right\}$ be the unit disk with the boundary $\gamma=\left\{\zeta\in\mathbb{C},\ |\zeta|=1\right\}$, and as a positive direction we choose the counter-clockwise direction, as usual. Let $S$ be a finite domain in the complex $z$ plane bounded by a simple smooth closed contour $L$, and let
\begin{equation*}
\omega : \zeta\in\mathbb{D}\mapsto z=\omega(\zeta)\in S
\end{equation*}
be a mapping, which maps $S$ onto $\mathbb{D}$ in the plane $\zeta$. The function $\omega(\zeta)$ is a holomorphic function inside of $\gamma$.\par
By taking complex conjugation of the second equation in Kolosov-Muskhelishvili formulae~(\ref{kolosov_formulae_under_mapping}), the following relation is obtained:
\begin{equation}
\label{normal_and_tangential_stresses_polar_coordinates_mapping}
\sigma_{rr}-i\,\sigma_{r\varphi} = \Phi(\zeta)+\overline{\Phi(\zeta)} - \frac{\zeta^{2}}{r^{2}\overline{\omega'(\zeta)}}\left[\overline{\omega(\zeta)}\Phi'(\zeta) + \omega'(\zeta)\,\Psi(\zeta)\right].
\end{equation}\par
For transforming the boundary value problem of linear elasticity in the unit disk of the $\zeta$-plane into a Riemann-Hilbert boundary value problem for a holomorphic function, discontinuities of holomorphic functions defined on $\mathbb{C}\setminus\gamma$ need to be described. Therefore, we consider the exterior of the unit disk $\mathbb{E}:=\mathbb{C}\setminus\mathbb{D}$, and we introduce holomorphic reflections as follows
\begin{equation}
l_{R} \colon \zeta\in\mathbb{E}\mapsto \overline{l\left(\frac{1}{\overline{\zeta}}\right)},
\end{equation}
where the index $R$ stays for reflection of a function and will be used in the sequel. The function $l_{R}(\zeta)$ is holomorphic in $\mathbb{E}$, if the function $l(\zeta)$ is holomorphic in $\mathbb{D}$.\par 
Let now $t=e^{i\varphi}\in\gamma$ be a point of the unit circle, and let $t^+$ and $t^-$ tend to $t$ from the interior and exterior of the unit disk, respectively. Thus, $t^+$ and $t^-$ can be defined as follows
\begin{equation}
\label{tplustminus}
\left\{
\begin{array}{l}
t^{+} := r^{+}e^{i\varphi}\in\mathbb{D} \mbox{ for } r^+\underset{r^{+}<1}{\rightarrow} 1, \\
\displaystyle {t^{-} := \frac{1}{\overline{t^{+}}}}\in\mathbb{E}.
\end{array}\right.
\end{equation}\par
Let now $\gamma_{\sigma}$ denotes the part of boundary $\gamma$, where traction boundary conditions are defined. Note that $\gamma_{\sigma}$ can be a union of several disjoint arcs, see \cite{Muskhelishvili_1,Muskhelishvili_2} for details. Considering relations $\overline{l(t^+)}=\overline{l(1/\overline{t^-})}=l_R(t^-)$, equation ~(\ref{normal_and_tangential_stresses_polar_coordinates_mapping}) can be now rewritten for a point $t\in\gamma_{\sigma}$ as follows
\begin{equation}
\label{eq:boundaryconditionstress}
\begin{array}{lcl}
\displaystyle \overline{\omega'(t)}\left[\sigma^*_{rr}-i\,\sigma^*_{r\varphi}\right] & = & \displaystyle \omega_R'(t^-)\Phi_R(t^-)+\omega_R'(t^-)\Phi(t^+) \\
\\
& & \displaystyle - t^2\left[\omega_R(t^-)\Phi'(t^+) + \omega'(t^+)\,\Psi(t^+)\right],
\end{array}
\end{equation}
where $\omega_R'$ is the reflection of the derivative function, and the left-hand side represents a known stress function on the boundary $\gamma_{\sigma}$ with $\sigma^*_{rr}$ and $\sigma^*_{r\varphi}$ being imposed stresses on $\gamma_{\sigma}$.\par 
To formulate a classical Riemann-Hilbert problem for a holomorphic function we need to rewrite equation~(\ref{eq:boundaryconditionstress}) in terms of only one holomorphic function, rather than a combination of several functions as it is written at the moment. For that we need to introduce an additional assumption: the conformal mapping has to be holomorphic on the entire complex plane $\mathbb{C}$ and not only on the unit disk $\mathbb{D}$
\begin{equation*}
\omega:\zeta\in\mathbb{C}\mapsto\omega(\zeta).
\end{equation*}
Consequently, $\omega_R(\zeta)$ is also defined on the entire complex plane and, in particular, in the interior of the unit disk $\mathbb{D}$. Hence, $\omega_R(t^-)$ and $\omega_R'(t^-)$ can be replaced by $\omega_R(t^+)$ and $\omega_R'(t^+)$ in~(\ref{eq:boundaryconditionstress}), respectively. It should be noted that if $\omega(\zeta)$ has a pole at infinity of order not higher than $N$, then the asymptotic expansion $\omega(\zeta)$ at the infinity can be written as follows
\begin{equation*}
\omega(\zeta)\underset{\left|z\right|\rightarrow +\infty}{=}\sum_{k=0}^N\omega_k\zeta^k.
\end{equation*}
In addition, if $\omega(\zeta)$ has a pole at infinity, then $\omega_R(\zeta)$ has a pole at zero and the asymptotic expansion have the form
\begin{equation}
\label{eq:pole}
\omega_R(\zeta)\underset{\left|z\right|\rightarrow 0}{=}\sum_{k=0}^N\frac{\overline{\omega}_k}{\zeta^k}.
\end{equation}\par
Let us consider the following holomorphic function on $\mathbb{C}\setminus\gamma$ 
\begin{equation}
\label{eq:omegadef}
\Omega:\zeta\in\mathbb{C}\setminus\gamma\mapsto\left\{
\begin{array}{ll}
\displaystyle{\omega_R'(\zeta)\Phi(\zeta)-\zeta^2\left[\omega_R(\zeta)\Phi'(\zeta) + \omega'(\zeta)\,\Psi(\zeta)\right]}, & \mbox{if } \left|\zeta\right|<1,\\
\\
\displaystyle -\omega_R'(\zeta)\Phi_R(\zeta), & \mbox{if } \left|\zeta\right|>1,
\end{array}
\right.
\end{equation}
where the origin has been removed from the domain for the case if $\omega_R(\zeta)$ has a pole at zero. However, if $\omega_R(\zeta)$ does not have a pole at the origin, then the origin should be added to the domain. The boundary condition (\ref{eq:boundaryconditionstress}) can now be written as
\begin{equation}
\overline{\omega'(t)}\left[\sigma^*_{rr}(t)-i\,\sigma^*_{r\varphi}(t)\right] =\Omega(t^+)-\Omega(t^-).
\end{equation}\par
Similar to $\gamma_{\sigma}$, we denote by $\gamma_{u}$ the part of $\gamma$, where displacements are prescribed. Again, $\gamma_{u}$ can be a union of several disjoint arcs. From~(\ref{displacement_cartesian}) and by help of variables $t^+$ formula for displacement boundary condition for a point $t\in\gamma_{u}$ can be written as follows
\begin{equation*}
2\mu(u^*_1-iu^*_2)=\kappa\overline{\eta(t^+)}-\overline{\omega(t^+)}\Phi(t^+)-\chi(t^+),
\end{equation*}
where $u^*_1$ and $u^*_2$ are known displacements along Cartesian directions in the original domain $S$ imposed on $\gamma_u$ considered as a function of $\varphi$. Finally, we need a formula for $\left(u_{1}^{*}\right)' - i\left(u_{2}^{*}\right)'$ with
\begin{equation*}
\left(u_{1}^{*}\right)' = \frac{\partial u_{1}^{*}}{\partial\varphi}, \qquad \left(u_{2}^{*}\right)' = \frac{\partial u_{2}^{*}}{\partial\varphi}.
\end{equation*}
Differentiating the previous formula we obtain
\begin{equation}
\label{eq:conditiondisplacementlast}
\begin{array}{lcl}
\displaystyle -2\mu i t\left[\left(u_{1}^{*}\right)'-i\left(u_{2}^{*}\right)'\right] & = & \displaystyle \overline{\omega'(t^+)}\Phi(t^+)- \kappa\overline{\omega'(t^+)}\overline{\Phi(t^+)} \\
\\ 
& & \displaystyle - t^2\left(\overline{\omega(t^+)}\Phi'(t^+)+\omega'(t^+)\Psi(t^+)\right),
\end{array}
\end{equation}
where the relations $\eta'(\zeta)=\omega'(\zeta)\Phi(\zeta)$ and $\chi'(\zeta)=\omega'(\zeta)\Psi(\zeta)$ have been used. Taking into account the assumption that $\omega(\zeta)$ is defined on $\mathbb{C}$, we finally get
\begin{equation*}
\begin{array}{lcl}
\displaystyle -2\mu i t\left[\left(u_{1}^{*}\right)'-i\left(u_{2}^{*}\right)'\right] & = & \displaystyle \omega_R'(t^+)\Phi(t^+)-\kappa\omega_R'(t^-)\Phi_R(t^-) \\
\\
& & \displaystyle - t^2\left(\omega_R(t^+)\Phi'(t^+)+\omega'(t^+)\Psi(t^+)\right),
\end{array}
\end{equation*}
or in terms of function~(\ref{eq:omegadef}),
\begin{equation}
-2\mu i t\left[\left(u_{1}^{*}\right)'-i\left(u_{2}^{*}\right)'\right] = \Omega(t^+)+\kappa\Omega(t^-).
\end{equation}
Thus, we get the following Riemann-Hilbert problem for the holomorphic function $\Omega$
\begin{equation}
\label{eq:hilbertproblem}
\left\{
\begin{array}{rcl}
\displaystyle \Omega(t^+)-\Omega(t^-)&  = & f(t) \mbox{ on } \gamma_{\sigma},\\
\displaystyle \Omega(t^+)+\kappa\Omega(t^-) & = & f(t) \mbox{ on } \gamma_u,
\end{array}
\right.
\end{equation}
with boundary function $f(t)$ defined by
\begin{equation}
\label{boundary_function}
f(t):=\left\{
\begin{array}{lcl}
\overline{\omega'(t)}\left[\sigma^*_{rr}-i\,\sigma^*_{r\varphi}\right]&\text{on}&\gamma_{\sigma},\\[3mm]
\displaystyle{-2\mu i t\left[\left(u_{1}^{*}\right)'-i\left(u_{2}^{*}\right)'\right]}&\text{on}&\gamma_u.
\end{array}
\right.
\end{equation}
It should be noted that the imposed normal and tangential stresses $\sigma^*_{rr}$ and $\sigma^*_{r\varphi}$ on $\gamma_{\sigma}$ correspond to polar directions in the $\zeta$-plane, although the imposed displacements $u^*_1$ and $u^*_2$ on $\gamma_u$ correspond to Cartesian directions in the $z$-plane. However, both imposed stresses and displacements are seen as functions of $\varphi$, or equivalently $t$, in the $\zeta$-plane.\par

\subsection{Solution of the Riemann-Hilbert boundary value problem in the unit disk for a general case}

In this section we describe at first the solution of Riemann-Hilbert problem~(\ref{eq:hilbertproblem}) for a general case, and later we specify it for the considered problem. Consider that $\gamma_{u}$ is the union of $n$ arcs such as $\gamma_{u}=\cup_{k=1}^n(a_k,b_k)$. Let us consider the following holomorphic function on $\mathbb{C}\setminus\gamma_u$
\begin{equation}
\label{particular_solutions}
X_{0}:\zeta\in\mathbb{C}\setminus\gamma_u\mapsto \prod\limits_{k=1}^{n}(\zeta-a_{k})^{-\frac{1}{2}+i\beta}(\zeta-b_{k})^{-\frac{1}{2}-i\beta}, \mbox{ with } \beta=\frac{\log\kappa}{2\pi}.
\end{equation}
Taking into account displacement and traction boundary conditions given on $\gamma_{u}$ and $\gamma_{\sigma}$, it is well known that the following classical relations hold \cite{Muskhelishvili_1}:
\begin{equation*}
\frac{X_0(t^+)}{X_0(t^-)}=-\kappa\text{ on }\gamma_{u} \mbox{ and } \frac{X_0(t^+)}{X_0(t^-)}=1\text{ on }\gamma_{\sigma}.
\end{equation*}
Thus, mixed boundary value problem~(\ref{eq:hilbertproblem}) can be reduced to the following problem for $\frac{\Omega(\zeta)}{X_0(\zeta)}$
\begin{equation}
\label{eq:cauchyproblem}
\frac{\Omega(t^+)}{X_0(t^+)}-\frac{\Omega(t^-)}{X_0(t^-)}=\frac{f(t)}{X_0(t^+)}\mbox{ on  }\gamma=\gamma_{\sigma}\cup\gamma_{u}.
\end{equation}
Solution of~(\ref{eq:cauchyproblem}) requires describing asymptotic behaviour of $\frac{\Omega(\zeta)}{X_0(\zeta)}$. For that, we recall that function $\Phi(\zeta)$ is holomorphic in $\mathbb{D}$, and therefore, we have
\begin{equation*}
\left\{\begin{array}{rcl}
\displaystyle \Phi(\zeta) & = & \displaystyle \sum_{k=0}^{+\infty}A_k\zeta^k, \mbox{ if } \left|\zeta\right|<1,\\
\\
\displaystyle\Phi_R(\zeta) & = & \displaystyle \sum_{k=0}^{+\infty}\frac{\overline{A}_k}{\zeta^k}, \mbox{ if } \left|\zeta\right|>1,\\
\end{array}\right.
\end{equation*}
where $A_{k},k=0,1,\ldots$ are unknown coefficients of the decomposition. Next, using the fact that $\omega'_R(\zeta)\underset{\left|\zeta\right|\rightarrow +\infty}{\rightarrow}\omega'(0)$ and definition of $\Omega(\zeta)$, we obtain the following asymptotic expansion
\begin{equation*}
\Omega(\zeta)\underset{\left|\zeta\right|\rightarrow +\infty}{=} B_{0}+\frac{B_{1}}{\zeta}+\frac{B_{2}}{\zeta^2}+\cdots,
\end{equation*}
where $B_{0},B_{1},\ldots$ are unknown coefficients. The asymptotic expansion of $1/X_0(\zeta)$ is obtained from~(\ref{particular_solutions}) as follows
\begin{equation}
\frac{1}{X_0(\zeta)}\underset{\left|\zeta\right|\rightarrow +\infty}{=}\zeta^n+D_{n-1}\zeta^{n-1}+\cdots+D_{1}\zeta+D_0+\frac{D_{-1}}{\zeta}+\cdots,
\end{equation}
where $D_{n-1},\ldots,D_{0},\ldots$ are known coefficients obtained by an 
asymptotic expansion of $\frac{1}{X_0(\zeta)}$. Finally, it is evident that there exists a polynomial of degree not higher than $n$
\begin{equation}
P_{n}(\zeta)=\widehat{C}_{0} + C_{1}\zeta + \ldots + C_{n}\zeta^{n},
\end{equation}
such that
\begin{equation}
\label{eq:limits}
\frac{\Omega(\zeta)}{X_0(\zeta)}-P_{n}(\zeta)\underset{\left|\zeta\right|\rightarrow +\infty}{\rightarrow}0.
\end{equation}\par
If $\omega_R(\zeta)$ has a pole at the origin, then the asymptotic expansion of $\frac{\Omega(\zeta)}{X_0(\zeta)}$ at the origin has to be determined. If the order of this pole of $\omega_R(\zeta)$ is not higher than $N$, as it has been shown in~(\ref{eq:pole}), then the pole of $\omega_R'(\zeta)$ at the origin is not higher than $N-1$. Considering that the value of $X_{0}(\zeta)$ at the origin is a non-zero constant, it follows from~(\ref{eq:omegadef}) that there exists a function $Q_N(\zeta)$ of the form
\begin{equation*}
Q_N(\zeta)=\widetilde{C}_0+\frac{C_{-1}}{\zeta}+\cdots +\frac{C_{-(N-1)}}{\zeta^{N-1}},
\end{equation*}
such that
\begin{equation}
\label{eq:limitszero}
\frac{\Omega(\zeta)}{X_0(\zeta)}-Q_{N}(\zeta)\underset{\left|\zeta\right|\rightarrow 0}{\rightarrow}0.
\end{equation}
Introducing a new constant $C_0:=\widehat{C}_0+\widetilde{C}_0$, we finally obtain:
\begin{equation*}
R(\zeta)=P_n(\zeta)+Q_N(\zeta)=C_n\zeta^n+\cdots+C_0+\frac{C_{-1}}{\zeta}+\cdots+\frac{C_{-(N-1)}}{\zeta^{N-1}}.
\end{equation*}
Thus, the general solution of~(\ref{eq:hilbertproblem}) is given now by
\begin{equation}
\label{R-H_general_solution}
\Omega(\zeta)=\Omega_0(\zeta)+X_0(\zeta)R(\zeta) \mbox{ with } \Omega_0(\zeta)=\frac{X_0(\zeta)}{2i\pi}\int_{\gamma}\frac{f(t)\text{d}t}{X_0(t^+)(t-\zeta)},
\end{equation}
where the integration is taken over the whole boundary $\gamma$. The coefficients $C_{-(N-1)},\cdots,C_0,\cdots,C_n$ should be identified by ensuring displacement continuity at ends of the arcs $a_k$ and $b_k$ and by ensuring that there is no stress and displacement singularities at zero.\par
Finally, holomorphic functions $\Phi(\zeta)$ and $\Psi(\zeta)$ can be easily derived from (\ref{R-H_general_solution}) and therefore displacements and stresses are obtained in $S$.\par

\subsection{Solution of the Riemann-Hilbert boundary value problem in the unit disk for the considered crack configuration}

\begin{figure}[h!]
\centering
\includegraphics[width=1.0\textwidth]{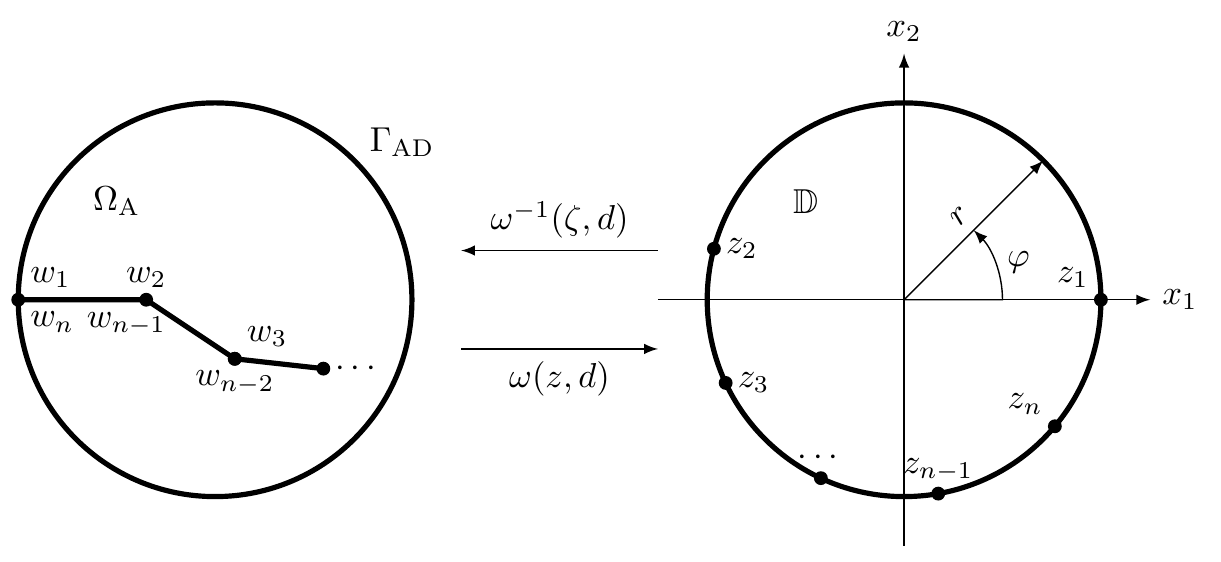}
\caption{Vertices and pre-vertices for the mapping between the unit disk and a cracked disk during the crack propagation process}
\label{figure_6}
\end{figure}

Next, we discuss the construction of an explicit solution of the Riemann-Hilbert problem for the crack propagation process shown in Fig. \ref{figure_3_crack_development}. Domain $\Omega_{\mathrm{A}}$ with a kinked crack  can be considered as a circular-arc polygon with vertices $w_{i}$, $i=1,\ldots n$, which are located along the crack path, and we keep the convention that vertex $w_{\frac{n+1}{2}}$ is located at the crack tip. Since according to the coupling idea, displacements are interpolated on the whole coupling interface $\Gamma_{\mathrm{AD}}$, no extra vertices are required on $\Gamma_{\mathrm{AD}}$, and the fact of having several coupling elements will be addressed in a piecewise definition of boundary function $f(t)$ in~(\ref{eq:hilbertproblem}). Thus, vertices $w_{i}$, $i=1,\ldots n$ are mapped to the corresponding pre-vertices at the unit circle $\gamma$ denoted by $z_{i}$, $i=1,\ldots n$, see Fig.~\ref{figure_6}.\par
Thus, in the case of analytical-numerical coupling the unit circle $\gamma$ is subdivided into arcs $\gamma_{u}=z_{n}z_{1}$ with unknown displacement boundary conditions given by the interpolation function~(\ref{interpolation_function_general_form}), and $\gamma_{\sigma} = \cup_{i=1}^{n-1} z_{i}z_{i+1}$ with traction-free conditions on the crack faces. Therefore, considering that only one arc with displacement boundary conditions is given and assuming that $\omega_R(\zeta)$ has no pole at the origin, a general solution~(\ref{eq:hilbertproblem}) can be written as follows
\begin{equation}
\label{general_solution_our_case}
\Omega(\zeta) = \frac{X_{0}(\zeta)}{2\pi i\kappa}\int\limits_{\gamma}\frac{f(t)dt}{X_{0}^{+}(t)(t-\zeta)} + X_{0}(\zeta)\left[C_{0}+C_{1}\zeta\right],
\end{equation}
where
\begin{equation*}
X_{0}(\zeta)=(\zeta-z_{1})^{-\frac{1}{2}-i\beta}(\zeta-z_{n})^{-\frac{1}{2}+i\beta}, \mbox{ with } \beta=\frac{\log\kappa}{2\pi}.
\end{equation*}
Applying displacement boundary conditions on $\Gamma_{\mathrm{AD}}$ and traction-free conditions on the crack faces, the following system of equations for unknown coefficients is obtained
\begin{equation}
\label{equations_for_unknown_coefficients}
\left\{\begin{array}{rcl}
\displaystyle (\kappa+1)\int\limits_{z_{n}z_{1}}\Omega_{0}(t_{0})\omega'(t_{0})dt_{0} \\
\\
\displaystyle + (\kappa+1)\int\limits_{z_{n}z_{1}}\omega'(t_{0})X_{0}(t_{0})(C_{0} + C_{1}t_{0})dt_{0} & = & \displaystyle 2\mu[f(z_{1})-f(z_{n})], \\
\\
\displaystyle \frac{X_{0}(0)}{2\pi\,i\kappa}\int\limits_{\gamma}\frac{f(t)}{X_{0}^{+}(t)}\frac{dt}{t} + C_{0}X_{0}(0) + \overline{C}_{1} & = & 0,
\end{array}\right.
\end{equation}
with
\begin{equation*}
\Omega_{0}(t_{0}) = \frac{X_{0}(t_{0})}{2\pi\,i\kappa}\int\limits_{\gamma}\frac{f(t)dt}{X_{0}^{+}(t)(t-t_{0})},
\end{equation*}
and
\begin{equation*}
\displaystyle X_{0}(0) = \lim\limits_{\zeta\to 0}\left[(\zeta-z_{1})^{-\frac{1}{2}-i\beta}(\zeta-z_{n})^{-\frac{1}{2}+i\beta}\right] =  e^{-i\pi}e^{-i\frac{\varphi_{1}+\varphi_{n}}{2}}e^{-\beta(\varphi_{n}-\varphi_{1})},
\end{equation*}
where the fact that $\ln|z_{1}|$ and $\ln|z_{n}|$ are zero on the unit disk has been taken into account. Denoting by $\varphi_{0}$ the argument of the middle of the arc $z_{n}z_{1}$, and by $\omega_{0}$ its central angle, the expression for $X_{0}(0)$ can be simplified to
\begin{equation*}
X_{0}(0)=-e^{-i\varphi_{0}-\beta\omega_{0}}.
\end{equation*}\par
To identify constants $C_{0}$ and $C_{1}$, system~(\ref{equations_for_unknown_coefficients}) can be transformed into its real form and solved explicitly. To avoid bulky expressions, we omit the presentation of the explicit solution of the corresponding real 4 by 4 system here. Nonetheless, the whole procedure remains the same on each step of the crack propagation process as long as boundary conditions of the Riemann-Hilbert problem kept as described in this section. The main computational complexity is related to numerical calculation of the conformal mapping.\par

\section{Energetic approach to crack propagation}\label{section_energetic_approach}

In this section, the mechanical point of view on the crack propagation based on the energetic formulation proposed in \cite{francfort-marigo} is described. Similar to previous sections, we describe a general setting of the energetic approach at first, and after that, we specify it for the problem considered in the paper.\par
Consider now time-dependent Neumann boundary conditions $\mathbf{F}(t)$ given on $\Gamma_{1}$, then for any time $t$, the crack geometry $\Gamma_{c}(t)$ is obtained by finding the global minimum of the total energy $\mathcal{E}^{\text{tot}}$ under the assumption of irreversibility of the crack growth, i.e. the crack can only grow. The total energy for Neumann boundary conditions $\mathbf{F}^{*}$ on $\Gamma_{1}$ and for any crack geometry $\Gamma_{c}^{*}$ is expressed as follows
\begin{equation*}
\mathcal{E}^{\text{tot}}(\mathbf{F}^{*},\Gamma_{c}^{*}) := \mathcal{E}(\mathbf{F}^{*},\Gamma_{c}^{*})-\mathcal{W}(\mathbf{F}^{*},\Gamma_{c}^{*})+\mathcal{D}(\Gamma_{c}^{*}),
\end{equation*}
where the stored elastic energy $\mathcal{E}(\mathbf{F}^{*},\Gamma_{c}^{*})$, the work of external forces $\mathcal{W}(\mathbf{F}^{*},\Gamma_{c}^{*})$, and the dissipated energy $\mathcal{D}(\Gamma_{c}^{*})$, are given by
\begin{equation*}
\begin{array}{lcl}
\displaystyle \mathcal{E}(\mathbf{F}^{*},\Gamma_{c}^{*}) & = & \displaystyle \frac{1}{2}\int_{G}\tensor{\sigma}(\mathbf{F}^{*},\Gamma_{c}^{*}):\tensor{\varepsilon}(\mathbf{F}^{*},\Gamma_{c}^{*})\text{d}V, \\
\\
\displaystyle \mathcal{W}(\mathbf{F}^{*},\Gamma_{c}^{*}) & = & \displaystyle \int_{\Gamma_{1}}\mathbf{F}^{*}\cdot\mathbf{u}(\mathbf{F}^{*},\Gamma_{c}^{*})\text{d}S, \quad \mathcal{D}(\Gamma_{c}^{*})=G_{c}\int_{\Gamma_{c}^{*}}\text{d}l,
\end{array}
\end{equation*}
where $\tensor{\sigma}$ is the stress tensor, $\tensor{\varepsilon}$ the strain tensor and $\mathbf{u}$ the displacement vector. Based on the above consideration, the energetic criterion can now be formulated as follows \cite{francfort-marigo}:
\begin{equation}
\label{eq:minimization}
\left\{\begin{array}{ll}
(a): & \displaystyle \forall\, s<t,\,\Gamma_c(s)\subset\Gamma_{c}(t),\\
\\
(b): & \displaystyle \forall\, \Gamma_{c}(t)\subset \Gamma_{c}^{*},\,\mathcal{E}^{\text{tot}}(\mathbf{F}(t),\Gamma_{c}(t))\leq \mathcal{E}^{\text{tot}}(\mathbf{F}(t),\Gamma_{c}^{*}), \\
\\
(c): & \displaystyle{\forall\, s<t,\,\mathcal{E}^{\text{tot}}(\mathbf{F}(t),\Gamma_{c}(t))\leq \mathcal{E}^{\text{tot}}(\mathbf{F}(t),\Gamma_{c}(s))}.
\end{array}\right.
\end{equation}
Let us make some remarks regarding the criterion: condition $(a)$ corresponds to the constraint of irreversibility of the crack growth; the condition $(b)$ ensures that the total energy for the actual crack is lower than for any longer crack; and condition $(c)$ ensures that the total energy of the actual crack is lower than for any previous real crack considering the actual boundary conditions.\par
Energetic criterion~(\ref{eq:minimization}) is formulated for a continuous time, in practice, however, a time discretisation $t_{1} < \cdots <t_{n}$ is introduced with $t_{1}$ corresponding to the initial configuration. Thus, according to~(\ref{eq:minimization}), knowing the crack geometry $\Gamma_{c}(t_{j-1})$ at the time step $j-1$, the crack geometry $\Gamma_{c}(t_j)$ at the time step $j$ (with $1\leq j\leq n$) is determined as follows 
\begin{equation}
\label{eq:minimization2}
\Gamma_c(t_j)=\underset{\Gamma_{c}(t_{j-1})\subset\Gamma_{c}^{*}}{\text{argmin}}\left[\mathcal{E}^{\text{tot}}(\mathbf{F}(t_{j}),\Gamma_{c}^{*})\right]
\end{equation}
Indeed, for any time discretisation, \eqref{eq:minimization2} clearly implies \eqref{eq:minimization}.\par
In general, the total energy $\mathcal{E}^{\text{tot}}(\mathbf{F}(t_{j}),\Gamma_{c}^{*})$ depends on the stress and displacement field in the whole domain $\Omega$. However, under the assumption that the crack can propagate only inside the analytical domain $\Omega_{\mathrm{A}}$, minimisation problem~(\ref{eq:minimization2}) can be formulated locally. In this case of local formulation, Neumann boundary conditions $\mathbf{F}_{\mathrm{A}}$ on the coupling interface $\Gamma_{\mathrm{AD}}$ are considered. These Neumann boundary conditions are obtained on each step of propagation $j$ and for each trial of new crack geometry by solving the continuous coupling with finite element method. Indeed, as the crack growth tends to relax strain and stress, the Neumann boundary conditions $\mathbf{F}_{\mathrm{A}}$ needs to be re-computed for any tested evolution of the crack geometry. However, as the analytical domain $\Omega_{\mathrm{A}}$ is chosen to cover the largest possible area in the elastic body, the computational cost is expected to be reduced. Such a local formulation would not be possible in the classical finite element setting without using elements of higher regularity, since traces of generalised derivatives of basis functions are needed in order to obtain Neumann data on $\Gamma_{\mathrm{AD}}$. However, this problem does not appear in the case of analytical-numerical coupling described in previous sections, since a strong solution to the differential equation in $\Omega_{\mathrm{A}}$ is constructed. Thus, Neumann data on $\Gamma_{\mathrm{AD}}$ can be obtained straightforwardly.\par
So, to formulate the minimisation problem, we consider Neumann boundary conditions $\mathbf{F}_{\mathrm{A}}$, which are formally determined as a function of $\mathbf{F}(t)$ and $\Gamma_{c}^{*}$, on the coupling interface $\Gamma_{\mathrm{AD}}$ for any crack $\Gamma_{c}^{*}$. Then, minimisation problem (\ref{eq:minimization2}) can be reduced to:
\begin{equation}
\label{eq:minimizationA}
\Gamma_c(t_j)=\underset{\Gamma_{c}(t_{j-1})\subset\Gamma_{c}^{*}}{\text{argmin}}\left[\mathcal{E}_{\mathrm{A}}^{\text{tot}}(\mathbf{F}_{\mathrm{A}}\left[\mathbf{F}(t_{j}),\Gamma_{c}^{*}\right],\Gamma_{c}^{*})\right],
\end{equation}
where {\itshape the local total energy} $\mathcal{E}_{\mathrm{A}}^{\text{tot}}(\mathbf{F}^{*}_{\mathrm{A}},\Gamma_{c}^{*})$ is given by
\begin{equation*}
\mathcal{E}_{\mathrm{A}}^{\text{tot}}(\mathbf{F}^{*}_{\mathrm{A}},\Gamma_{c}^{*})=\mathcal{E}_{\mathrm{A}}(\mathbf{F}^{*}_{\mathrm{A}},\Gamma_{c}^{*})-\mathcal{W}_{\mathrm{A}}(\mathbf{F}^{*}_{\mathrm{A}},\Gamma_{c}^{*})+\mathcal{D}(\Gamma_{c}^{*}),
\end{equation*}
with the elastic energy $\mathcal{E}_{\mathrm{A}}(\mathbf{F}^{*}_{\mathrm{A}},\Gamma_{c}^{*})$ stored in the analytical domain $\Omega_{\mathrm{A}}$, and the work of forces $\mathcal{W}_{\mathrm{A}}(\mathbf{F}^{*}_{\mathrm{A}},\Gamma_{c}^{*})$ on the coupling interface $\Gamma_{\mathrm{AD}}$ are given by
\begin{equation}
\label{eq:local_quantities}
\begin{array}{rcl}
\displaystyle \mathcal{E}_{\mathrm{A}}(\mathbf{F}^{*}_{\mathrm{A}},\Gamma_{c}^{*}) & = & \displaystyle \frac{1}{2}\int_{\Omega_{\mathrm{A}}}\tensor{\sigma}(\mathbf{F}^{*}_{\mathrm{A}},\Gamma_{c}^{*}):\tensor{\varepsilon}(\mathbf{F}^{*}_{\mathrm{A}},\Gamma_{c}^{*})\text{d}V, \\
\\
\displaystyle \mathcal{W}_{\mathrm{A}}(\mathbf{F}^{*}_{\mathrm{A}},\Gamma_{c}^{*}) & = & \displaystyle \int_{\Gamma_{\mathrm{AD}}}\mathbf{F}^{*}_{\mathrm{A}}\cdot\vect{u}(\mathbf{F}^{*}_{\mathrm{A}},\Gamma_{c}^{*})\text{d}S.
\end{array}
\end{equation}
Thus, formulation~\eqref{eq:minimizationA} presents the advantage that the analytical solution of the Riemann-Hilbert problem described in Section~\ref{section_conformal_mapping} is used to compute at each time step $j$ the total energy $\mathcal{E}_{\mathrm{A}}^{\text{tot}}$. Indeed, for all $\mathbf{F}^{*}_{\mathrm{A}}$ and $\Gamma_{c}^{*}$ one can compute analytically $\tensor{\sigma}(\mathbf{F}^{*}_{\mathrm{A}},\Gamma_{c}^{*})$, $\tensor{\varepsilon}(\mathbf{F}^{*}_{\mathrm{A}},\Gamma_{c}^{*})$ and  $\vect{u}(\mathbf{F}^{*}_{\mathrm{A}},\Gamma_{c}^{*})$ involved in~(\ref{eq:local_quantities}).\par

\section{First example towards a complete numerical scheme}\label{Section_conclusions}

The aim of this section is two-fold: at first, we briefly discuss the difficulties related to practical implementation of the complete solution strategy presented in this paper, and recall some of possible approaches to overcome these difficulties, which will constitute the future work; after that, we present a small numerical example focusing only on the use of conformal mapping and Riemann-Hilbert problem, since these are the crucial parts of the complete algorithm to model crack propagation in elastic bodies. Moreover, we underline openly all problems related to the numerical stability of the method, since overcoming these problems constitute the major part of future work.\par
The analytical domain $\Omega_{\mathrm{A}}$ containing a crack is a circular-arc polygon, with crack-faces representing the polygonal part and the coupling interface $\Gamma_{\mathrm{AD}}$ being the circular arc. The idea of the method presented in this paper is to map the circular-arc polygon to the unit disk, because Riemann-Hilbert problems in the unit disk are well studied in the context of linear elasticity, see for example \cite{Gakhov,Muskhelishvili_1}. While there are several classical works studying conformal mappings of circular arc-polygon regions, see for example \cite{Bjorstad,Driscoll} and references therein, it is well-known that an explicit representation of a mapping function between a circular-arc polygon and the unit disk does not exist. The classical approach to construct a mapping function for such type of domains is to work with the Schwarz-Christoffel differential equation.\par 
Because the Schwarz-Christoffel differential equation is ill-posed due to nonlinear constrains for the parameters of the map, it is known that its numerical solution is a challenging task, although some methods for numerical calculations of such mappings exist \cite{Andersson,Brown,Howell}. An alternative approach would be to use directly algorithms for numerical conformal mapping, such as for example the osculation algorithms \cite{Henrici,Porter}. However, the main obstacle for the use of numerical conformal mapping in the context of coupled method is the fact, that not only the geometry must be mapped, as typically addressed in the field of numerical conformal mappings, but the differential equation and its solution procedure as well. Thus, it must be studied how the solution of Riemann-Hilbert problem in our case will behave under numerical conformal mapping.\par
Because of difficulties discussed above on the way of implementing the complete numerical procedure presented in this paper, we present an illustrative example focusing only on the crack propagation based on the solution of Riemann-Hilbert problem. Thus, instead of considering a global boundary value problem in a domain $G$, we consider a boundary value problem formulated directly in the analytical domain $\Omega_{\mathrm{A}}$ and boundary conditions on the coupling interface $\Gamma_{\mathrm{AD}}$. Additionally, to avoid a circular-arc polygon mapping, we consider a square domain centred at the crack tip of the initial configuration.\par
\begin{figure}[h!]
\centering
\includegraphics[width=0.65\textwidth]{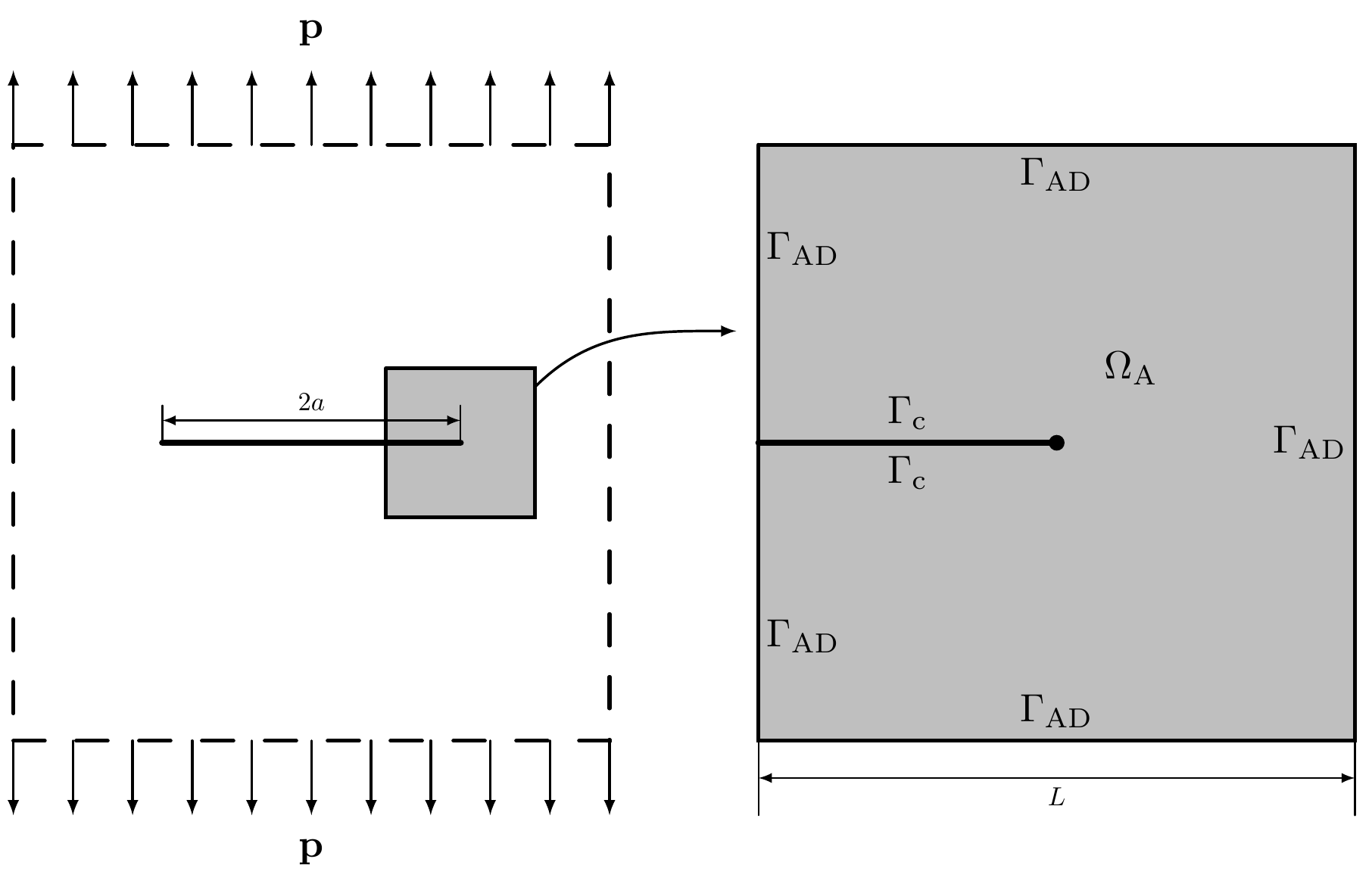}
\caption{Setting for the illustrative example: crack in an infinite body (left), representation of the analytical domain $\Omega_{\mathrm{A}}$ with the coupling interface $\Gamma_{\mathrm{AD}}$ as a square (right)}
\label{figure_test_example}
\end{figure}

Let us consider an infinite plane containing a single crack of a length $2a$ with constant stresses $\mathbf{p}$ applied at infinity (Fig.~\ref{figure_test_example}, left). To formulate a boundary value problem, we consider a square domain of length $L$ located around one of the crack tips (Fig.~\ref{figure_test_example}, right) representing the analytical domain $\Omega_{\mathrm{A}}$. To keep the illustrative example closer to the setting discussed in Section~\ref{section_conformal_mapping}, displacement boundary conditions are considered on the interface $\Gamma_{\mathrm{AD}}$ and traction-free conditions on the crack faces $\Gamma_{\mathrm{c}}$. Thus, we consider the following boundary value problem
\begin{equation*}
\left\{
\begin{array}{rcll}
-\mu\,\Delta\mathbf{u}-(\lambda+\mu)\,\mathrm{grad}\,\mathrm{div}\,\mathbf{u} & = & 0, & \mbox{in } \Omega_{\mathrm{A}},\\
\mathbf{u} & = & u_{1}+iu_{2}, & \mbox{on } \Gamma_{\mathrm{AD}},\\
\displaystyle \boldsymbol{\sigma}\cdot \overline{n} & = & 0, & \mbox{on } \Gamma_{\mathrm{c}}, 
\end{array}\right.
\end{equation*}
where the displacements components $u_{1}$ and $u_{2}$ are chosen according to the well-known analytical solution, see for example \cite{Liebowitz}, and are given by the following formulae:
\begin{equation*}
\begin{array}{rcl}
\displaystyle u_{1} & = & \displaystyle \frac{p\sqrt{2ra}}{8\mu}\left[(2\kappa-1)\cos\left(\frac{\varphi}{2}\right) - \cos\left(\frac{3\varphi}{2}\right)\right], \\
\\
\displaystyle u_{2} & = & \displaystyle\frac{p\sqrt{2ra}}{8\mu}\left[(2\kappa+1)\sin\left(\frac{\varphi}{2}\right) - \sin\left(\frac{3\varphi}{2}\right)\right],
\end{array}
\end{equation*}
with $r$, $\varphi$ being polar coordinates with the coordinate origin located at the crack tip, and $\kappa$ and $\mu$ being material parameters.\par
\begin{figure}[h!]
\centering
\includegraphics[width=0.98\textwidth]{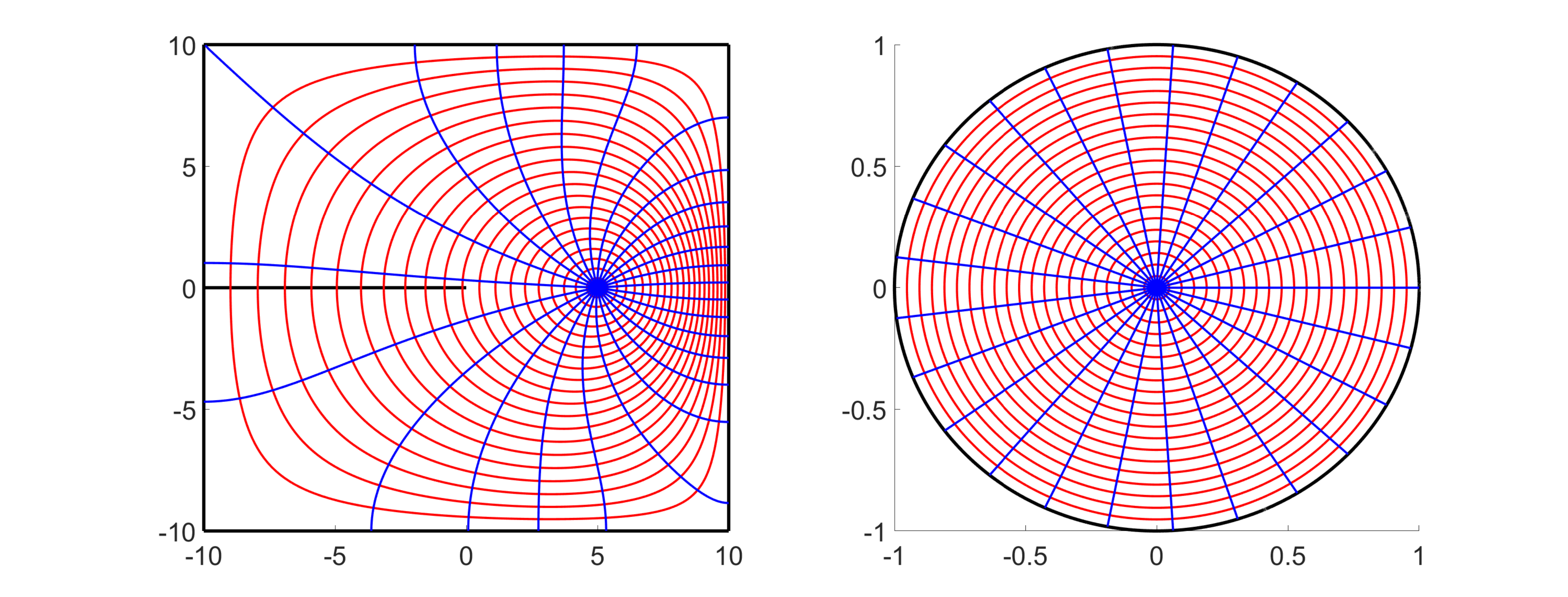}
\caption{Domain $\Omega_{\mathrm{A}}$ together with the conformal grid (left), and the unit disk with the conformal grid obtained after calculating the conformal mapping from $\Omega_{\mathrm{A}}$ by using the PlgCirMap toolbox}
\label{figure_conformal_map}
\end{figure}

For a numerical conformal mapping of the domain $\Omega_{\mathrm{A}}$, in general, the classical Schwarz-Christoffel toolbox for Matlab developed by T.A. Driscoll \cite{Driscoll_1} can be used. However, using this toolbox implies the necessity to work with the inverse Schwarz-Christoffel mapping in all constructions presented in Section~\ref{section_conformal_mapping}, which complicates the numerical part. Therefore, instead of the Schwarz-Christoffel toolbox, the PlgCirMap Matlab toolbox will be used, which has been introduced recently in \cite{Nasser}. The PlgCirMap toolbox allows mapping of polygonal multiply connected domains onto circular domains by using Koebe's iterative method. Fig.~\ref{figure_conformal_map} shows the domain $\Omega_{\mathrm{A}}$ and the unit disk together with the conformal grid calculated by the PlgCirMap toolbox.\par
The advantage of using the PlgCirMap toolbox is the fact that the direct mapping from the polygonal domain $\Omega_{\mathrm{A}}$ to the unit disk can be used in all constructions presented in Section~\ref{section_conformal_mapping}, which significantly simplifies all related calculations. Nonetheless, although the PlgCirMap toolbox provides a lot of useful functions for numerical conformal mapping, it is also not free of geometrical restrictions: polygonal domains with slits and cusps are not allowed. To overcome this restriction, we model the crack in a domain as a cut with width of order $10^{-4}$. In this case, the conformal mapping to the unit disk can be calculated with the relative residual of order $10^{-5}$. The vertices calculated during the conformal mapping, as well as pre-vertices of the original domain, are listed in Table~\ref{table_vertices}. Because some vertices are located very close to each other, we list the coordinates of vertices in a long format provided by Matlab.\par
\begin{table}[h!]
\begin{tabular}{|c|c|c|}
\hline
{\bfseries Number} & {\bfseries Pre-vertices} & {\bfseries Vertices} \\
\hline
1 & $-10 + 0.0001i$ & $0.955417215003076 + 0.295259115482937i$ \\
2 & $0.0001i$ & $0.955422389606861 + 0.295242370668431i$ \\
3 & $-0.0001i$ & $0.955391244897614 + 0.295343138015746i$ \\
4 & $-10 - 0.0001i$ & $0.955396421837903 + 0.295326390861588i$ \\
5 & $-10 - 10i$ & $0.825604379631623 + 0.564249420321443i$ \\
6 & $10 - 10i$ & $-0.171634834630707 + 0.985160638444964i$ \\
7 & $10 + 10i$ & $0.414173847997492 - 0.910197793688246i$ \\
8 & $-10 + 10i$ & $1$ \\
\hline
\end{tabular}
\caption{Vertices and pre-vertices for the conformal mapping}\label{table_vertices}
\end{table}\par
\begin{figure}[h!]
\centering
\includegraphics[width=0.98\textwidth]{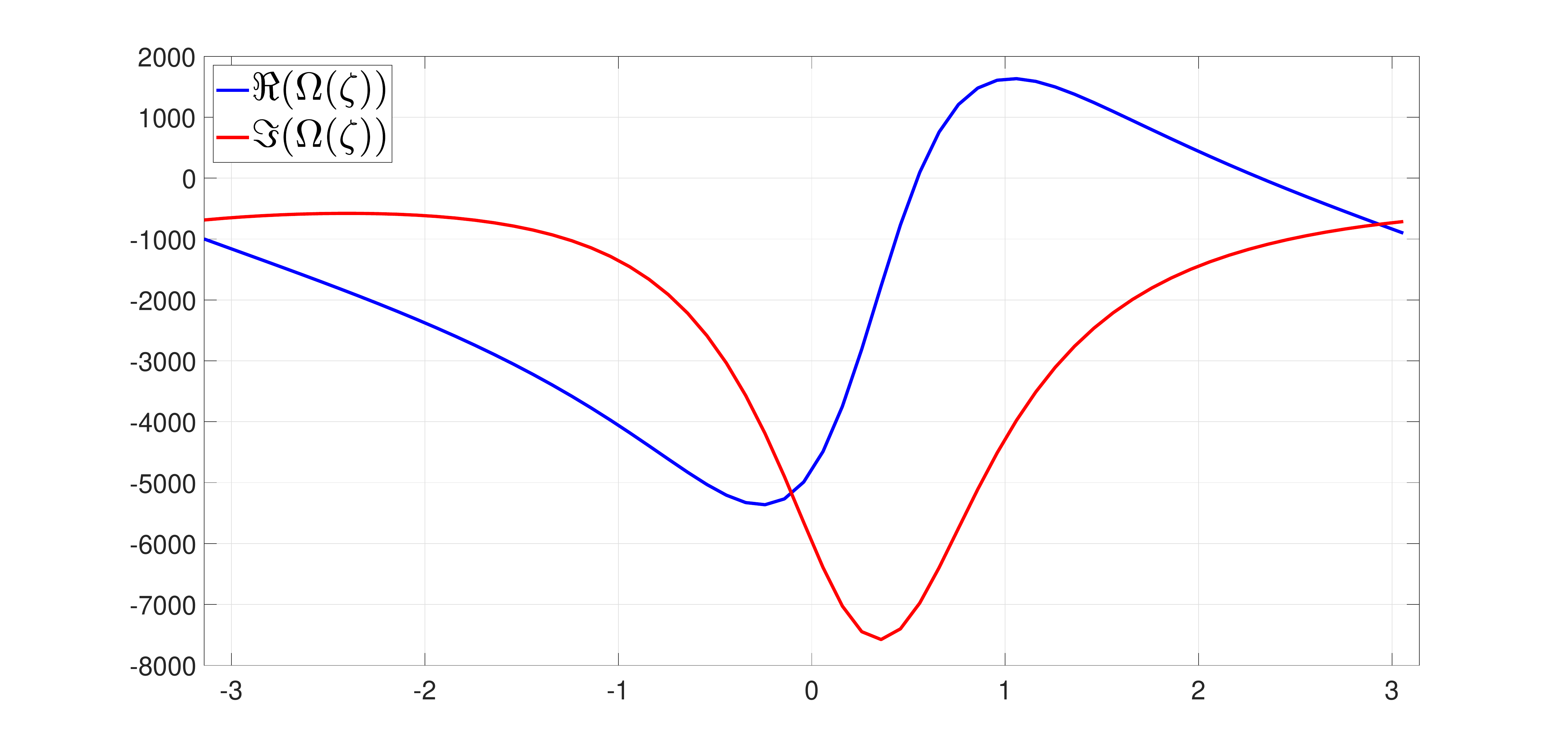}
\caption{Domain $\Omega_{\mathrm{A}}$ together with the conformal grid (left), and the unit disk with the conformal grid obtained after calculating the conformal mapping from $\Omega_{\mathrm{A}}$ by using the PlgCirMap toolbox}
\label{figure_omega_plot}
\end{figure}

Let us now outline the general procedure for constructing a solution of the Riemann-Hilbert problem:
\begin{itemize}
\item[Step 1.] Map the domain $\Omega_{\mathrm{A}}$ to the unit disk.
\item[Step 2.] Map boundary conditions from $\Omega_{\mathrm{A}}$ to the unit disk.
\item[Step 3.] Create and solve linear system of equations~(\ref{equations_for_unknown_coefficients}).
\item[Step 4.] Compute the general solution of Riemann-Hilbert problem in the unit disk by help of formula~(\ref{general_solution_our_case}).
\end{itemize}
Fig.~\ref{figure_omega_plot} shows the solution of Riemann-Hilbert problem in the unit disk with respect to $\varphi\in[-\pi,\pi]$ and for $r=\frac{1}{2}$. It is also important to remark, that the solution of a linear system on Step 3 can be written explicitly in our case, implying that no numerical procedure is necessary to solve the linear system. Nonetheless, computing the solution is still numerically difficult, because several singular integrals need to be calculated on Step 3, since they appear in the coefficients of the system and in the right-hand side. Thus, the quality of the solution of Riemann-Hilbert problem (and further computations with it) strongly depends on calculation of these singular integrals. However, because four of the eight vertices are located very close to each other, see Table~\ref{table_vertices}, they cause numerical stability issues during computing the singular integrals. In the example presented in this section, the singular integrals could be computed only with the accuracy of order $10^{-1}$ by using built-in Matlab adaptive quadratures. Evidently, this accuracy is not sufficient for further calculations of stresses and displacements. Therefore, one of the tasks for future work is finding a numerical quadrature for computing singular integrals with the accuracy of the same order as provided by numerical conformal mapping.\par

\section{Summary and outlook}\label{Section_summary}

In this paper, we have presented the theoretical background of a coupled analytical-numerical approach to model a crack propagation process in two-dimensional bounded domains. The main idea of the method is to obtain the correct solution behaviour near the crack tip by help of the analytical solution constructed by using tools of the complex function theory and couple it continuously with the finite element solution in the region far from singularity. To calculate possible directions of crack growth, the idea is to utilise the conformal mapping techniques and to transform a problem of linear elasticity into a Riemann-Hilbert problem in the unit disk for holomorphic functions. In the paper, we have presented the analytical solution of the Riemann-Hilbert problem, as well as discussed numerical stability issues appearing on the way of practical realisation of the method, proposed in this paper.\par
As it has been discussed in Section~\ref{Section_conclusions}, the main difficulty of the method is related to the need of having a conformal mapping between a circular-arc polygon and the unit disk. Unfortunately, this mapping cannot be expressed explicitly by help of known conformal mappings. Therefore, we have considered a simplified version of a problem in Section~\ref{Section_conclusions}, where a circular domain has been replaced by a rectangular domain. Nonetheless, even in that case, further studies are necessary for finding a numerical quadrature enabling calculating of singular integrals with a higher accuracy, which is necessary for calculating stresses and displacements.\par
In summary, this paper presents a work-in-progress, rather than a final result. The scope of future work consists in studying of different approaches for practical calculations of circular-arc polygon mappings in the context of coupling method, as well as analysing of different advanced methods for computing singular integrals. Additionally, further theoretical studies of the method, such as for example unique solvability of the interpolation problem under conformal mapping, must be also made.\par
Finally, it is worth to mention, that some recent works dealing with analysis of kinked cracks proposed to work with a mapping from a half-space \cite{Adda_Bedia,Beom}. Considering that different conformal mappings can be used on different propagation steps, as well as a composition of several conformal mappings can also be helpful in practice, the use of the mapping from a half-space needs also to be studied in the context of coupled method, presented in this paper. Perhaps a new setting for the complete methods can be found on this way.\par


\begin{thebibliography}{99}

\bibitem{Adda_Bedia}
M. Adda-Bedia, R. Arias, \emph{Brittle fracture dynamics with arbitrary paths I. Kinking of a dynamic crack in general antiplane loading}. Journal of the Mechanics and Physics of Solids, 51, pp. 1287-1304, 2003.

\bibitem{Andersson}
A. Andersson, \emph{A modified Schwarz-Christoffel mapping for regions with piecewise smooth boundaries}. Journal of Computational and Applied Mathematics, 213, pp. 56-70, 2008.

\bibitem{Beom}
H.G. Beom, J.W. Lee, C.B. Cui, \emph{Analysis of a kinked crack in an anisotropic material under antiplane deformation}. Journal of Mechanical Science and Technology, 26(2), pp. 411-419, 2012.

\bibitem{Bjorstad}
P. Bj\o rstad, E. Grosse, \emph{Conformal mapping of circular arc polygons}. SIAM Journal on Scientific and Statistical Computing, 8(1), pp. 19-32, 1987.

\bibitem{Brown}
P.R. Brown, \emph{Mapping onto circular arc polygons}, Complex Variables. 

Theory and Application: An International Journal, 50:2, pp. 131-154, 2005.

\bibitem{Ciarlet}
P.G. Ciarlet, \emph{The finite element method for elliptic problems}. North-Holland, Amsterdam, 1978.

\bibitem{Driscoll_1}
T.A. Driscoll, \emph{A Matlab toolbox for Scwarz-Christoffel mapping}. ACM Transactions on Mathematical Software, 22, pp. 168-186, 1996. 

\bibitem{Driscoll}
T.A. Driscoll, L.N. Trefethen, \emph{Scwarz-Christoffel mapping}. Cambridge University Press 2002.

\bibitem{erdogan1963}
F. Erdogan, G.C. Sih, \emph{On the crack extension in plates under plane loading and transverse shear}. Journal of basic engineering, 85(4), pp. 519-525, 1963. 

\bibitem{francfort-marigo}
G.A. Francfort, J.-J. Marigo, \emph{Revisiting brittle fracture as an energy minimization problem}. Journal of the Mechanics and Physics of Solids, 46(8), pp. 1319–1342, 1998.

\bibitem{Gakhov}
F.D. Gakhov, \emph{Boundary value problems}. Pergamon Press, 1966.

\bibitem{Guerlebeck_1}
K. G\"urlebeck, D. Legatiuk, \emph{On the continuous coupling of finite elements with holomorphic basis functions}. Hypercomplex Analysis: New perspectives and applications, ISBN 978-3-319-08770-2, Birkh\"auser, Basel, 2014.

\bibitem{Guerlebeck_2}
K. G\"urlebeck, U. K\"ahler, D. Legatiuk, \emph{Interpolation problem arising in a coupling of finite element method with holomorphic basis functions}. AIP Conference proceedings, Volume 1648, 2015. DOI: 10.1063/1.4912655.

\bibitem{Guerlebeck_3}
K. G\"urlebeck, U. K\"ahler, D. Legatiuk, \emph{Error estimates for the coupling of analytical and numerical solutions}. Complex Analysis and Operator Theory, 11(5), pp. 1221-1240, 2017.

\bibitem{Henrici}
P. Henrici, \emph{A general theory of osculation algorithms for conformal mapping}. Linear Algebra and its Applications, 52/53, pp. 361-382, 1983.

\bibitem{Howell}
L.H. Howell, \emph{Numerical conformal mapping of circular arc polygons}. Journal of Computational and Applied Mathematics, 46, pp. 7-28, 1993.

\bibitem{leblond1989crack}
J.B. Leblond, \emph{Crack paths in plane situations—I. General form of the expansion of the stress intensity factors}. International Journal of Solids and Structures, 25(11), pp. 1311-1325, 1989.

\bibitem{Legatiuk_1}
D. Legatiuk, H.M. Nguyen, \emph{Improved convergence results for the finite element method with holomorphic functions}. Advances in Applied Clifford Algebra, 24(4), pp. 1077-1092, 2014.

\bibitem{Legatiuk_2}
D. Legatiuk, \emph{Evaluation of the coupling between an analytical and a numerical solution for boundary value problems with singularities}. PhD Thesis, Bauhaus-Universit\"at Weimar, 2015. ISBN: 978-3-95773-193-7

\bibitem{Liebowitz}
H. Liebowitz, \emph{Fracture, an advanced treatise. Volume II: Mathematical fundamentals}. Academic Press, 1968.

\bibitem{Moes99}
N. Mo{\"e}s, J. Dolbow, T. Belytschko, \emph{A finite element method for crack growth without remeshing}. International journal for numerical methods in engineering, 46(1), pp. 131-150, 1999.

\bibitem{Muskhelishvili_1}
N.I. Muskhelishvili, \emph{Singular integral equations: boundary problems of functions theory and their applications to mathematical physics}. Wolters-Noordhoff Publishing, 1958.

\bibitem{Muskhelishvili_2}
N.I. Muskhelishvili, \emph{Some basic problems of the mathematical theory of elasticity}. Springer Science+Business Media Dordrecht, 1977.

\bibitem{Nasser}
M.M.S. Nasser, \emph{PlgCirMap: A MATLAB toolbox for computing conformal mappings from polygonal multiply connected domains onto circular domains}. SoftwareX 11, 100464, 2020.

\bibitem{Piltner_1}
R. Piltner, \emph{Special finite elements with holes and internal cracks}. International Journal of Numerical Methods in Engineering, 21, pp. 509-528, 1985.

\bibitem{Piltner_2}
R. Piltner, \emph{The derivation of special purpose element functions using complex solution representation}. Computer Assisted Mechanics and Engineering Sciences, 10(4), pp. 597-607, 2003.

\bibitem{Porter}
R.M. Porter, \emph{An accelerated osculation method and its application to numerical conformal mapping}. Complex Variables, Theory and Application: An International Journal, 48, pp. 569-582, 2003.

\bibitem{sif}
I.N. Sneddon, M. Lowengrub, \emph{Crack problems in the classical theory of elasticity}. John Willey \& Sons, 1969.

\bibitem{MAXG}
C.H. Wu, \emph{Fracture under combined loads by maximum-energy-release-rate criterion}. Journal of Applied Mechanics, 45(9), pp. 553-558, 1978.
 
\end{thebibliography}
\end{document}